%Paper: dg-ga/9412003
%From: huebschm@mpim-bonn.mpg.de (Johannes Huebschmann)
%Date: Wed, 21 Dec 1994 22:44:31 --100

%&amstex
%This paper is written in AMSTeX
\documentstyle{amsppt}
\magnification=1200

\hoffset=-0.5pc
\vsize=57.2truepc
\hsize=38truepc
\nologo
\spaceskip=.5em plus.25em minus.20em

 \define\atiboton{1}

\define\goldmone{6}

\define\narastwo{20}
\define\narashed{21}
\define\sjamlerm{22}
\define\weinstwe{23}

\define\sbaueron{24}
\define\baueroko{25}
\define\biereck{26}
\define\bodenone{27}
\define\kbrownon{28}
\define\donalthr{29}
\define\furustee{30}
\define\guhujewe{31}
\define\modus{32}
\define\aspheric{33}
\define\einsther{34}
\define\huebjeff{35}
\define\kirkklas{36}
\define\magnusbo{37}
\define\mehtrama{38}
\define\mehtsesh{39}
\define\orvogzie{40}
\define\patteone{41}
\define\seshaone{42}
\define\seshaboo{43}
\define\weiltwo{44}
\noindent
dg-ga/9412003
\newline
\noindent
\topmatter
\title
Symplectic and Poisson structures
of certain moduli spaces. II. \\
Projective representations of
cocompact
planar discrete groups
\endtitle
\author Johannes Huebschmann{\dag}
\endauthor
\affil
Max Planck Institut f\"ur Mathematik
\\
Gottfried Claren Str. 26
\\
D-53 225 BONN
\\
huebschm\@mpim-bonn.mpg.de
\endaffil
\abstract{
Let $G$ be a Lie group with a biinvariant metric, not necessarily
positive definite. It is shown that a certain construction carried out
in an earlier paper for the fundamental group of a closed surface may be
extended to an arbitrary infinite orientation preserving  cocompact planar
discrete group of euclidean or non-euclidean motions $\pi$ and yields (i) a
symplectic structure on a certain smooth manifold $\Cal M$ containing the
space $\roman{Hom}(\pi,G)$ of homomorphisms and, furthermore, (ii) a
hamiltonian  $G$-action on $\Cal M$ preserving the symplectic structure
together with a momentum mapping in such a way that the reduced space
equals the space $\roman{Rep}(\pi,G)$ of representations. More generally,
the construction also applies to certain spaces of projective
representations. For $G$ compact, the resulting spaces of representations
inherit structures of {\it stratified symplectic space\/} in such a way that
the strata have finite symplectic volume . For example,
{\smc Mehta-Seshadri} moduli spaces of semistable holomorphic parabolic
bundles with rational weights or spaces closely related to them
arise in this way  by {\it symplectic reduction in finite dimensions\/}.}
\endabstract
%\date{February 4, 1994} \enddate
\thanks{\dag} The author carried out this work in the framework of
 the VBAC research group of EUROPROJ.
\endthanks
\endkeywords
\endsubjclass

\endtopmatter
\document
\leftheadtext{Johannes Huebschmann}
\rightheadtext{Projective representations of cocompact planar discrete groups}
\beginsection Introduction

Let  $\pi$  be a finitely generated
orientation preserving infinite  cocompact planar discrete group
of euclidean or hyperbolic motions.
So $\pi$
acts by isometries on the euclidean or
upper half plane (as appropriate)
and the orbit space $\Sigma$
is a compact orientable Riemann surface.
In the hyperbolic case, $\pi$ is also called a (cocompact)
{\it Fuchsian group\/}.
Let
$G$ be a Lie group with a biinvariant metric,
not necessarily positive definite.
The set
of representations
or more generally
projective representations
of $\pi$ in $G$
is known to inherit additional structure,
under suitable circumstances,
cf. e.~g.
\cite\goldmone,
\cite\narastwo,
\cite\narashed,
\cite\sbaueron,
\cite\baueroko,
\cite\bodenone,
\cite\furustee,
\cite\modus,
\cite\huebjeff,
\cite\kirkklas,
\cite\mehtsesh,
\cite\seshaone,
\cite\seshaboo,
\cite\weiltwo.
In this paper
we shall
study
the
symplectic or more generally {\it Poisson\/} geometry
of such representations.
More specifically,
extending a certain construction
carried out in an earlier paper \cite\modus\
for the fundamental group
of a closed surface,
we shall obtain certain {\it extended representation spaces\/}
(see below for a precise definition).
In \cite \modus\
the theory has been made for arbitrary
finite presentations
but it has been applied only to the standard
presentation of the fundamental group of a closed surface.
However, the general approach  in
\cite \modus\
applies to an arbitrary
group of the kind $\pi$
and yields the following, see (2.9) below for a more
precise statement.

\proclaim{Theorem}
There is a smooth symplectic
manifold $\Cal M$
containing
$\roman{Hom}(\pi,G)$
(as a deformation retract)
together with
a hamiltonian  $G$-action
on $\Cal M$
and momentum mapping
in such a way that the reduced space
equals
the space $\roman{Rep}(\pi,G)$
of representations of $\pi$ in $G$.
More generally,
spaces of projective representations of $\pi$ are obtained
by symplectic reduction at appropriate non-zero values of the momentum
mapping.
\endproclaim

The (smooth) symplectic manifold
$\Cal M$
together with the $G$-action and momentum mapping
is what we mean by an {\it extended representation space\/}.
The second clause of the Theorem will be made precise in Section 3 below.
Our result, apart from being interesting in its own right,
reveals some interesting and attractive geometric properties of these
twisted representation spaces,
which have been spelled out
in \cite\modus\
for the special case
considered there.
For example,
it implies that, for $G$ compact,
the resulting twisted representation spaces
inherit a structure of stratified symplectic space.
In particular it entails
that {\it symplectic reduction, applied to
a certain smooth finite dimensional symplectic
manifold with a hamiltonian action of
the unitary group\/} yields
certain stratified symplectic spaces
containing
as top stratum
the stable part of
the {\smc Mehta-Seshadri}
\cite\mehtsesh\
moduli spaces of semistable
parabolic vector
bundles with rational weights.
In fact the two spaces
are presumably globally homeomorphic
 but details have not been worked out yet.
For parabolic degree zero,
they {\it are\/} known to be homeomorphic,
by a result of {\smc Mehta-Seshadri}, cf.
(4.1) and (4.3) in \cite\mehtsesh.
\smallskip
The
projective representations
mentioned above, also referred to as {\it twisted representions\/}
below,
are
representations of certain
central extensions of $\pi$.
These extensions include
in particular
all fundamental groups of Seifert fibered  spaces
with empty boundary
(as 3-manifolds)
which are Eilenberg-Mac Lane spaces and have
orientable decomposition surface,
cf. \cite\orvogzie;
see Section 3 below for details.
Because of their relevance for Floer homology,
$\roman{SU}(N)$-representation spaces of Seifert fibered homology 3-spheres
have been studied in
\cite\sbaueron,
\cite\baueroko,
\cite\bodenone,
and \cite\kirkklas.
Among others, our approach
compactifies
spaces of irreducible representations
of fundamental groups
of the Seifert fibered spaces
belonging to the above class,
including the homology 3-spheres just mentioned,
as symplectic
manifolds
to  stratified symplectic spaces.
\smallskip
The  group $\pi$
is determined by its genus $\ell$ and
its torsion numbers $m_1,\dots,m_n$.
The basic idea is that the construction
carried out in our paper \cite\modus\
for a surface group
still works for such a  group $\pi$
when the free group on the generators is replaced
by a free product of a free group
of rank $2\ell$
together with
cyclic groups
$Q_1,\dots,Q_n$ of finite orders
respectively
$m_1,\dots,m_n$.
The key observation is
that
$\pi$, being assumed infinite
but cocompact,
is a two-dimensional Poincar\'e duality group
over the reals
\cite\biereck.
In particular,
the second homology group
$\roman H_2(\pi,\bold R)$
is a one-dimensional real vector space
and,
starting from a generator,
we can carry
out a construction of closed
2-form on an extended moduli space
similar to that in \cite\modus\
for the case of an ordinary surface group.
The infinitesimal structure
is then handled by means
of a corresponding length two {\it projective\/}
resolution
of the reals $\bold R$ in the category of
left $\bold R\pi$-modules,
and the fact that $\pi$
is a Poincar\'e duality group
over the reals
entails the nondegeneracy of the resulting
2-form, much as in \cite\modus.
\smallskip
I am indebted to A. Weinstein for discussions.
Any unexplained notation is the same as that
in our paper \cite\modus.

\beginsection 1. Poincar\'e duality for cocompact planar groups

The  group $\pi$
is given by a presentation
$$
\Cal P = \langle x_1,y_1,  \dots,x_\ell,y_\ell,z_1,\dots,z_n; r,
r_1,\dots,r_n\rangle,
$$
where
$$
r = \Pi [x_j,y_j] z_1 \dots z_n,
\ r_j = z_j^{m_j},\ m_j\geq 2.
$$
Since $\pi$ is assumed cocompact there are no parabolic elements
in the Fuchsian case.
The hypothesis
that $\pi$ is infinite is equivalent to the requirement that
the {\it measure\/}
$$
\mu(\pi) = 2\ell-2+\sum\left(1 - \frac 1{m_j}\right)
$$
be non-negative.
For example the group with the smallest positive measure is
given by the presentation
$$
\langle z_1,z_2,z_3; z_1z_2z_3, z_1^2, z_2^3, z_3^7 \rangle.
$$
All this is classical and may be found in standard textbooks,
see e.~g.~\cite\magnusbo.
\smallskip
The cohomology of $\pi$ is well known,
cf.  \cite \aspheric,
\cite \patteone, \cite\weiltwo.
In fact,
for an arbitrary ground ring $R$,
the Fox calculus, applied to the presentation $\Cal P$,
yields
a free resolution
$$
\bold R(\Cal P)\colon
\dots
@>>>
R_j(\Cal P)
@>\partial_j>>
\dots
@>>>
R_2(\Cal P)
@>\partial_2>>
R_1(\Cal P)
@>\partial_1>>
R\pi
$$
of $R$ in the category of left
$R\pi$-modules,
cf. our paper
\cite \aspheric.
We recall that, for $ j \geq 3$,
$R_j(\Cal P) =
R\pi[r_1,\dots,r_n]$,
the free left
$R\pi$-module
on the relators $r_1,\dots,r_n$,
while
$R_2(\Cal P) =R\pi[r,r_1,\dots,r_n]$
and
$R_1(\Cal P) =
R\pi[x_1,y_1,\dots,x_\ell,y_\ell,z_1,\dots,z_n]$.
Moreover the boundary operators
$\partial_j$ are given by the Fox calculus; in particular
$\partial_2$
is given by the matrix
of Fox derivatives,
which roughly amounts to reexpressing the elements
$r-1,r_1-1, \dots, r_n -1$
of the group ring $R\pi$ in terms of
$x_1-1, y_1-1, z_1-1$ etc.,
and
$\partial_1$ maps the generators $x_j$ etc. to
the elements $x_j-1$ etc. of $R\pi$.
In particular, the chain complex calculating the homology of $\pi$
with values in $R$ looks like
$$
\dots @>>> R[r,r_1,\dots,r_n]
@>\overline \partial_2>>
R[x_1,y_1,\dots,x_\ell,y_\ell,z_1,\dots,z_n]
@>{0}>>
R,
$$
where
the operator $\overline \partial_2$ satisfies the formulas
$$
\align
\overline \partial_2 (r) &= z_1 + \dots + z_n
\\
\overline \partial_2 (r_j) &= m_jz_j,
\quad 1 \leq j \leq n.
\endalign
$$
\smallskip
Henceforth we denote by $m$  the least common multiple of $m_1,\dots,m_n$.
The 2-chain
$$
b =mr - \frac m{m_1} r_1 - \dots - \frac m{m_n} r_n
$$
is then a 2-cycle,
and a closer look at the resolution reveals that
$\roman H_2(\pi,\bold Z)$
is then infinite cyclic, generated
by the class of $b$.
Likewise,
$\roman H_2(\pi,\bold R)$
is
a one-dimensional real vector space;
we shall take
$\kappa  = \frac 1m [b]  \in \roman H_2(\pi,\bold R)$
as its generator.
Moreover, the group $\pi$ is in fact a
{\it two-dimensional Poincar\'e duality
group\/} over the reals
having fundamental class
$\kappa \in \roman H_2(\pi,\bold R)$,
that is,
for every $\bold R \pi$-module $A$,
cap product with $\kappa$
yields a natural isomorphism
$$
\cap \kappa
\colon
\roman H^*(\pi,A)
@>>>
\roman H_{2-*}(\pi,A),
$$
cf. \cite\biereck.
In particular,
$\roman H^2(\pi,\bold R)$
is also a one-dimensional real vector space.
The group $\pi$ is thus of type $FP_2$ over the reals,
cf. e.g. \cite\kbrownon.
Moreover it has
$\roman H^1(\pi,\bold R\pi)$
zero.
In fact, $\pi$ has a torsion free
subgroup of finite index and hence
a
torsion free
normal subgroup of finite index,
necessarily a surface group, say $\tau$.
The cohomology
$\roman H^1(\pi,\bold R\pi)$
then amounts to the invariants
$\roman H^0(Q,\roman H^1(\tau,\bold R\pi))$,
with reference to the quotient group $Q = \pi \big / \tau$
which is finite.
However, as a $\tau$-module,
$\bold R\pi$ decomposes as a finite direct sum of copies
of
$\bold R\tau$,
and
$\roman H^1(\tau,\bold R\tau)$
is zero since
$\tau$ is a surface group.
Hence
$\roman H^1(\pi,\bold R\pi)$
is zero.
It is well known that a finitely presented
group
of type $FP_2$
having
$\roman H^1(\pi,\bold R\pi)$
zero
is a two-dimensional duality group.
The above considerations show that
the dualizing module is in fact that of the reals,
with trivial
$\pi$-module structure.
Hence
$\pi$
is a two-dimensional
Poincar\'e duality group
over the reals.
\smallskip
For
a ring $R$ containing
the rationals,
the free resolution
$\bold R(\Cal P)$
projects onto a projective resolution of length 2.
For later reference, we spell it out  for $R = \bold R$, the reals.
Let $j = 1,\dots, n$, and write
$Q_j$ for the finite cyclic subgroup
of $\pi$ generated by $z_j$;
it has exact order $m_j$, see e.~g. \cite\aspheric.
Since $Q_j$ is a finite group,
its augmentation ideal
$IQ_j \subseteq \bold RQ_j$
is a projective
$\bold RQ_j$-module
whence
the induced $\bold R\pi$-module
$\bold R\pi \otimes_{\bold RQ_j} IQ_j $
is projective.
We now consider
the beginning
$$
\bold R\pi[r,r_1,\dots,r_n]
@>\partial_2>>
\bold R\pi[x_1,y_1,\dots,x_\ell,y_\ell,z_1,\dots,z_n]
@>\partial_1>>
\bold R\pi
$$
of our free resolution
$\bold R(\Cal P)$
of the reals
and divide out
the generators $r_1,\dots,r_n$
of
$\bold R\pi[r,r_1,\dots,r_n]$
and their $\partial_2$-images
in
$\bold R\pi[x_1,y_1,\dots,x_\ell,y_\ell,z_1,\dots,z_n]$.
This yields a projective resolution
$$
\bold P(\Cal P)
\colon
P_2(\Cal P)
@>\partial_2>>
P_1(\Cal P)
@>\partial_1>>
\bold R\pi
$$
of $\bold R$
in the category of
left $\bold R\pi$-modules
with
$P_2(\Cal P)= \bold R\pi[r]$ and
$$
P_1(\Cal P) =
\bold R\pi[x_1,y_1,\dots,x_\ell,y_\ell]
\oplus
\bold R\pi \otimes_{\bold RQ_1} IQ_1
\oplus
\dots
\oplus
\bold R\pi \otimes_{\bold RQ_n} IQ_n.
$$
Notice, for $n=0$,
the group $\pi$ is just a surface group in the usual sense
and $\bold P(\Cal P)$ boils down to the
usual free resolution for a surface group.
\smallskip
For illustration, when $\ell = 1$ and $n = 1$, we have $r=[x,y]z$ and
$$
\partial_2 [r]
=
(1-xyx^{-1})[x]
+
(x-[x,y])[y]
+
[x,y][z].
$$

\beginsection 2. Representation spaces of planar discontinuous groups

Write $F$ for the free group
on the generators in $\Cal P$.
As in \cite\modus,
for a group $\Pi$ we denote by $(C_*(\Pi,R),\partial)$
the chain complex of its inhomogeneous reduced normalized
bar resolution over $R$.
Pick  $c\in C_2(F,\bold R)$
so that
$\partial c =  [r]- \frac 1{m_1} [r_1] - \dots - \frac 1{m_n} [r_n]
\in C_1(F,\bold R)$;
this is certainly possible.
Modulo the normal closure $N$ of the relators,
$c$ will then represent the class $\kappa$.
As in \cite\modus,
let $O$ be the subspace of the Lie algebra $g$ where the exponential
mapping is regular.
Denote by $\rho$ the word map from
$G^{2\ell+n} =\roman{Hom}(F,G)$
to
$G^{1+n}$
for the presentation $\Cal P$,
and
write
$\Cal H(\Cal P,G)$
for the
smooth manifold determined by the requirement that
a pull back diagram
$$
\CD
\Cal H(\Cal P,G)
@>\widehat \rho>> O  \times O^n
\\
@V{\eta}VV
@VV{\roman{exp}}V
\\
\roman{Hom}(F,G)
@>>\rho> G \times G^n
\endCD
$$
of spaces results, where the induced map
from
$\Cal H(\Cal P,G)$
to
$O^{1+n}$ is denoted by $\widehat \rho$.
The construction
in our paper \cite\modus\
yields
(i) a 2-chain $c \in C_2(F,\bold R)$ (\cite\modus\ Lemma 2),
(ii) the equivariant
closed 2-form $\omega_{c,\Cal P} = \eta^*(\omega_c) - \widehat \rho^*B$
on $\Cal H(\Cal P,G)$ (\cite\modus\ Theorem~1), and
(iii) a smooth equivariant map
$\mu \colon \Cal H(\Cal P,G) \to g^*$
whose adjoint
$\mu^\sharp$ from $g$ to
$C^{\infty}(\Cal H(\Cal P,G))$ satisfies the identity
$$
\delta_G(\omega_{c,\Cal P}) = d \mu^\sharp
\colon
g @>>> \Omega^1(\Cal H(\Cal P,G))
$$
on $\Cal H(\Cal P,G)$ (\cite\modus\ Theorem 2).
In fact, we can at first carry out these constructions
with reference to $b$ and thereafter divide by $m$.
In particular,
$\omega_{c,\Cal P} - \mu^\sharp$
is an equivariantly closed form in
$(\Omega_G^{*,*}(\Cal H(\Cal P,G));d,\delta_G)$
of total degree 2.
Thus, cf. \cite\atiboton\ and
what is said in Section 5 of \cite \modus,
$\mu$
formally satisfies the property of being
a momentum mapping for
the 2-form $\omega_{c,\Cal P}$ on
$\Cal H(\Cal P,G)$, with reference to the
obvious $G$-action, except
that
$\omega_{c,\Cal P}$
is not necessarily non-degenerate.
When the standard homotopy operator
on forms on $g$ is taken,
the relevant map $\psi$ from $g$ to
$g^*$ is in fact the adjoint of the chosen 2-form
$\cdot$ on $g$,
cf.
the remark in Section 1 of \cite\modus.
The map $\mu$ then amounts to the composite of
$\widehat \rho$ with
the sum map from $O^{n+1}$ to $g$,
combined with the adjoint
of the given 2-form $\cdot$ on $g$.
\smallskip
Let $\phi \in \roman{Hom}(F,G)$,
and suppose that
$\phi(r)$ and
each
$\phi(r_j)$ lie in the centre of $G$.
Then the composite of $\phi$
with the adjoint representation of $G$
induces a structure of a $\pi$-module on $g$,
and we write
$g_{\phi}$ for $g$, viewed as as $\pi$-module in
this way.
Recall that the
2-form $\cdot$ on $g$
and the homology class $\kappa$
determine the alternating
2-form
$$
\omega_{\kappa,\cdot,\phi}
\colon
\roman H^1(\pi,g_{\phi})
\otimes
\roman H^1(\pi,g_{\phi})
@>{\cup}>>
\roman H^2(\pi,\bold R)
@>{\cap \kappa}>> \bold R
\tag2.1
$$
on
$\roman H^1(\pi,g_{\phi})$.
Application of the functor
$\roman{Hom}_{\bold R\pi}(\cdot\,, g_{\phi})$
to the free resolution
$\bold R(\Cal P)$
yields the chain complex
$$
\bold C(\Cal P, g_{\phi})\colon
\roman C^0(\Cal P, g_{\phi})
@>{\delta_{\phi}^0}>>
\roman C^1(\Cal P, g_{\phi})
@>{\delta_{\phi}^1}>>
\roman C^2(\Cal P, g_{\phi}),
\tag2.2
$$
cf. \cite\modus\ (4.1),
computing the group cohomology $\roman H^*(\pi,g_{\phi})$
in degrees 0 and 1; we recall that
there are canonical isomorphisms
$$
\roman C^0(\Cal P, g_{\phi}) \cong g,
\quad
\roman C^1(\Cal P, g_{\phi}) \cong g^{2\ell+ n},
\quad
\roman C^2(\Cal P, g_{\phi}) \cong g^{1+n} .
$$
To explain the geometric significance
of this chain complex,
denote by $\alpha_\phi$
the smooth map
from $G$  to $\roman{Hom}(F,G)$
which assigns $x \phi x^{-1}$ to $x \in G$,
and write
\linebreak
$R_\phi\colon g^{2\ell+n} \to  \roman T_{\phi} \roman{Hom}(F,G)$
and
$R_{\rho\phi}\colon g^{1+n} \to  \roman T_{\rho \phi}G^{1+n}$
for the corresponding operations of right translation.
The tangent maps
$\roman T_e\alpha_{\phi}$
and $\roman T_{\phi}\rho$
make commutative the diagram
$$
\CD
\roman T_eG
@>\roman T_e\alpha_{\phi}>>
\roman T_{\phi} \roman{Hom}(F,G)
@>{\roman T_{\phi} \rho}>>
\roman T_{\rho (\phi)}G^{1+n}
\\
@A{\roman{Id}}AA
@A{\roman R_{\phi}}AA
@A{\roman R_{\rho(\phi)}}AA
\\
g
@>>{\delta^0_{\phi}}>
g^{2\ell +n}
@>>{\delta^1_{\phi}}>
g^{1+n},
\endCD
\tag2.3
$$
cf. \cite\modus\ (4.2).
The commutativity of the diagram
(2.3) shows at once that
right translation
identifies the
kernel of
the derivative
$\roman T_{\phi} \rho$
with the kernel
of the coboundary operator
$\delta_{\phi}^1$
from
$\roman C^1(\Cal P, g_{\phi})$
to
$\roman C^2(\Cal P, g_{\phi})$,
that is, with the
vector space
$\roman Z^1(\pi,g_{\phi})$
of $g_{\phi}$-valued 1-cocycles of $\pi$;
this space does {\it not\/}
depend on a specific presentation $\Cal P$, whence the notation.
We note that $\roman C^1(\Cal P, g_{\phi}) =
Z^1(F,g_{\phi})$,
the space of $g_{\phi}$-valued 1-cocycles for $F$.
Pick $\widehat \phi \in \Cal H(\Cal P,G)$
so that the canonical map
$\eta$ from
$\Cal H(\Cal P,G)$
to
$G^{2\ell+n}$
sends
$\widehat \phi$ to $\phi$.
Then, likewise,
right translation identifies the
restriction of the
2-form $\omega_{c,\Cal P}$
to the kernel of
the derivative
$\roman T_{\widehat \phi} \widehat \rho$
with
the 2-form
on $\roman Z^1(\pi,g_{\phi})$
obtained as the
composite of
$\omega_{\kappa,\cdot,\phi}$
with the projection
from
$\roman Z^1(\pi,g_{\phi})$
to
$\roman H^1(\pi,g_{\phi})$,
cf.  Corollary 4.8 of \cite\modus.
\smallskip
To obtain the space $\roman{Hom}(\pi,G)$
as the zero locus of a momentum mapping
in the {\it usual\/} sense, we must
cut the space
$\roman{Hom}(F,G)$ to size, in the following way:
Denote by
$F^\natural$
the group
given by the presentation
$$
\langle x_1,y_1,  \dots,x_\ell,y_\ell,z_1,\dots,z_n;
z_1^{m_1},\dots,z_n^{m_n}\rangle.
$$
Its space
$
\roman{Hom}(F^\natural,G)
$
of homomorphisms
decomposes as
$$
\roman{Hom}(F^\natural,G)
\cong
G^{2\ell} \times
\roman{Hom}(Q_1,G)
\times
\dots
\times
\roman{Hom}(Q_n,G),
$$
and the generators
induce an embedding
of
$\roman{Hom}(F^\natural,G)$
into
$\roman{Hom}(F,G)$
as a smooth submanifold
with a finite number of connected components.
We shall say more about these connected components later.
By construction, the word map $\rho$
induces a smooth map
$r$ from
$\roman{Hom}(F^\natural,G)$
to
$G$
so that
$\roman{Hom}(\pi,G) = r^{-1}(e)$,
the pre-image of the neutral element $e$ of $G$.
This map merely depends on the relator $r$ whence
we denote it by that same symbol.
Application of the functor
$\roman{Hom}_{\bold R\pi}(\cdot, g_{\phi})$
to the projective resolution
$\bold P(\Cal P)$
then yields the chain complex
$$
\bold C(\bold P(\Cal P), g_{\phi})\colon
\roman C^0(\bold P(\Cal P), g_{\phi})
@>{\delta_{\phi}^0}>>
\roman C^1(\bold P(\Cal P), g_{\phi})
@>{\delta_{\phi}^1}>>
\roman C^2(\bold P(\Cal P), g_{\phi}),
\tag2.4
$$
computing the group cohomology $\roman H^*(\pi,g_{\phi})$
in all degrees; by construction,
there are canonical isomorphisms
$$
\roman C^0(\bold P(\Cal P), g_{\phi}) \cong g,
\quad
\roman C^2(\bold P(\Cal P), g_{\phi}) \cong g.
$$
Moreover, the canonical projection from
$\bold R(\Cal P)$ to
$\bold P(\Cal P)$
induces a canonical injection of
$\bold C(\bold P(\Cal P), g_{\phi})$
into
$\bold C(\Cal P, g_{\phi})$
identifying the former with a subcomplex
of the latter.
\smallskip
We now suppose that
$\phi$ lies in
$\roman{Hom}(F^\natural,G)$
viewed as a subspace
of
$\roman{Hom}(F,G)$, that is,
$\phi (r_j) $ is trivial for $j=1,\dots,n$ but, beware,
$\phi(r)$ is still admitted to be an arbitrary element of the centre of
$G$.
The following observation
will be crucial.

\proclaim{Proposition 2.5}
The tangent maps
$\roman T_e\alpha_{\phi}$
and $\roman T_{\phi}r$
make commutative the diagram
$$
\CD
\roman T_eG
@>\roman T_e\alpha_{\phi}>>
\roman T_{\phi} \roman{Hom}(F^\natural,G)
@>{\roman T_{\phi} r}>>
\roman T_{r (\phi)}G
\\
@A{\roman{Id}}AA
@A{\roman R_{\phi}}AA
@A{\roman R_{r(\phi)}}AA
\\
g
@>>{\delta^0_{\phi}}>
\roman C^1(\bold P(\Cal P), g_{\phi})
@>>{\delta^1_{\phi}}>
g,
\endCD
\tag2.6
$$
having its vertical arrows isomorphisms
of vector spaces.
\endproclaim

In fact, the diagram (2.3)
restricts to the diagram (2.6).

\smallskip
The next step is to carry out a construction
of a space similar to that of
the space denoted by $\Cal H(\Cal P,G)$ in \cite\modus.
To this end, write
$\Cal H^\natural(\Cal P,G)$
for the
space determined by the requirement that
a pull back diagram
$$
\CD
\Cal H^\natural(\Cal P,G)
@>\widehat r>> O
\\
@V{\eta}VV
@VV{\roman{exp}}V
\\
\roman{Hom}(F^\natural,G)
@>>r> G
\endCD
\tag2.7
$$
results, where the induced map
from $\Cal H^\natural(\Cal P,G)$ to $O$ is denoted by $\widehat r$.
The space
$\Cal H^\natural(\Cal P,G)$
is a smooth manifold
and the induced map
$\eta$ from
$\Cal H^\natural(\Cal P,G)$
to
$\roman{Hom}(F^\natural,G)$
is a smooth codimension zero immersion whence
$\Cal H^\natural(\Cal P,G)$ has the same dimension as
$\roman{Hom}(F^\natural,G)$;
moreover
the above injection
of $\roman{Hom}(\pi,G)$
into
$\roman{Hom}(F,G)$
induces a canonical
injection
of $\roman{Hom}(\pi,G)$
into
$\Cal H^\natural(\Cal P,G)$.
Further,
$\Cal H^\natural(\Cal P,G)$
may be viewed as a subspace
of $\Cal H(\Cal P,G)$ in a canonical way,
and, for the present $\phi$, the above chosen
$\widehat \phi \in \Cal H(\Cal P,G)$
actually lies in
$\Cal H^\natural(\Cal P,G)$.
The commutativity of the diagram
(2.6) then shows that
right translation
identifies the
kernel of
the derivative
$\roman T_{\phi} r$
with the
vector space
$\roman Z^1(\pi,g_{\phi})$
of $g_{\phi}$-valued 1-cocycles of $\pi$,
and the same is true of the
kernel of
the derivative
$$
\roman T_{\widehat \phi} \widehat r
\colon
\roman T_{\widehat \phi}
\Cal H^\natural(\Cal P,G)
@>>>
\roman T_{\widehat r \widehat \phi} O.
$$
Still by Corollary 4.8 of \cite\modus,
right translation identifies the
restriction of the
2-form $\omega_{c,\Cal P}$
to the kernel of
the derivative
$\roman T_{\widehat \phi} r$
with
the 2-form
on $\roman Z^1(\pi,g_{\phi})$
obtained as the
composite of
$\omega_{\kappa,\cdot,\phi}$
with the projection
from
$\roman Z^1(\pi,g_{\phi})$
to
$\roman H^1(\pi,g_{\phi})$.
\smallskip
Suppose that the
given 2-form $\cdot$ on $g$ is non-degenerate.
Write $Z$ for the centre of $G$ and $z$ for its Lie algebra.
By Poincar\'e duality
in the cohomology of $\pi$, for every
$\phi$
in the pre-image $r^{-1}(z) \subseteq \Cal H(\Cal P,G)$
of $z \subseteq O$, in particular,
for every $\phi \in \roman{Hom}(\pi,G)$,
the 2-form
$\omega_{\kappa,\cdot,\phi}$
is then symplectic.
Notice that when
$\pi$ is torsion free, that is,
a surface group,
this form boils down to that
considered by {\smc Goldman} \cite\goldmone\
and reexamined in \cite\modus.
\smallskip
Next we consider the restriction map
from
$(\Omega_G^{*,*}(\Cal H(\Cal P,G));d,\delta_G)$
to
\linebreak
$(\Omega_G^{*,*}(\Cal H^\natural(\Cal P,G));d,\delta_G)$
induced by the injection
of
$\Cal H^\natural(\Cal P,G)$
into
$\Cal H(\Cal P,G)$.
Abusing notation, we write
$\omega_{c,\Cal P}$
and
$\mu^{\sharp}$
for the classes
in
$(\Omega_G^{*,*}(\Cal H^\natural(\Cal P,G));d,\delta_G)$
obtained by restriction
to $(\Omega_G^{*,*}(\Cal H(\Cal P,G));d,\delta_G)$
of the corresponding classes denoted by the same symbol.
It is clear that
the identity
$
\delta_G (\omega_{c,\Cal P})
= d \mu^{\sharp}
$
passes to
$(\Omega_G^{*,*}(\Cal H^\natural(\Cal P,G));d,\delta_G)$.
In other words,
for every $X \in g$, we have
$$
\omega_{c,\Cal P}(X_{\Cal H^\natural(\Cal P,G)},\cdot\,)
= d (X \circ \mu),
\tag2.8
$$
that is, formally the momentum mapping property is satisfied,
perhaps up to a sign depending
on the choice of conventions
which is unimportant for the geometry
of the situation.
The rest of the construction is now exactly the same
as that in our paper \cite\modus, except that
we work with the space
$\Cal H^\natural(\Cal P,G)$
instead of $\Cal H(\Cal P,G)$
and
with
$\bold P(\Cal P)$
instead of
$\bold R(\Cal P)$.
In fact, the formal momentum mapping property (2.8),
together with
the symplecticity of the 2-form
$\omega_{\kappa,\cdot,\phi}$
at every
$\phi \in
r^{-1}(z)$,
implies that
$\omega_{c,\Cal P}$
has maximal rank
equal to
$\dim \Cal H^\natural(\Cal P,G)$
at every point
$\widehat \phi$ of
$\Cal H^\natural(\Cal P,G)$
in
the pre-image
$\widehat r^{-1}(z)$ of $z$, in particular, at every point of
$\roman{Hom}(\pi,G)$.
In fact,
given a point $\widehat \phi$
of $\widehat r^{-1}(z)$,
with the notation $\phi = \eta(\widehat \phi)$,
the symplecticity of the 2-form
$\omega_{\kappa,\cdot,\phi}$
implies that
the 2-form $\omega_{c,\Cal P}$
on
the tangent space
$\roman T_{\widehat \phi}\Cal H^\natural(\Cal P,G)
\cong\roman C^1(\bold P(\Cal P), g_{\phi})$,
restricted to the subspace
$\roman Z^1(\pi,g_{\phi})$ of 1-cocycles,
has degeneracy space
equal to the subspace $\roman B^1(\pi,g_{\phi})$ of 1-coboundaries,
and the momentum mapping property then implies that
the 2-form $\omega_{c,\Cal P}$
on
the whole space
$\roman C^1(\bold P(\Cal P), g_{\phi})$ is non-degenerate.
Let
$\Cal M^\natural(\Cal P,G)$
be the subspace of
$\Cal H^\natural(\Cal P,G)$
where
the 2-form $\omega_{c,\Cal P}$
is non-degenerate;
this is an open $G$-invariant
subset containing
the pre-image $\widehat r^{-1}(z)$.
Summarizing, we
obtain the following.

\proclaim{Theorem 2.9}
The space $\Cal M^\natural(\Cal P,G)$
is a smooth $G$-manifold
(having
in general more
than one connected component),
the 2-form $\omega_{c,\Cal P}$ is a $G$-invariant symplectic
structure on it,
and the restriction
$$
\mu = -\psi \circ r \colon
\Cal M^\natural(\Cal P,G)
@>>>
g^*
$$
is a momentum mapping
in the usual sense.
\endproclaim

Notice for a surface group the smooth manifold
$\Cal M^\natural(\Cal P,G)$
boils down
to the one denoted by
$\Cal M(\Cal P,G)$ in \cite\modus.
We already remarked that when the
constructions in \cite\modus\ are carried out
by means of the standard homotopy operator
on forms on $g$,
the  map $\psi$ from $g$ to
$g^*$ is  the adjoint of the chosen 2-form
$\cdot$ on $g$.
Whatever homotopy operator has been
taken, symplectic reduction then yields the space
$\roman{Rep}(\pi,G)$
of representations of $\pi$ in $G$.

\medskip\noindent{\bf 3. Twisted representation spaces}\smallskip\noindent
The construction of the
universal central extension
of the fundamental group of a closed surface
generalizes in the following way to our group $\pi$:
Let $F$ be the free group on the generators
of the presentation $\Cal P$,
$N$ the normal closure of $r, r_1,\dots, r_n$ in $F$,
and
$\Gamma = F\big / [F,N]$;
then the
kernel
$N\big /[F,N]$
 of the canonical projection from $\Gamma$ to
$\pi$
decomposes into a direct sum
of $n+1$ infinite cyclic groups,
generated by
$[r] = r[F,N]\in \Gamma$
and
$[r_1] = r_1[F,N],\dots,[r_n] = r_n[F,N]$;
this is an immediate consequence of the statement
of the {\it Identity Theorem\/} for the presentation $\Cal P$,
see e.~g. our paper \cite\aspheric.
Notice the decomposition
of $N\big /[F,N]$
depends on the presentation $\Cal P$, and {\it not\/}
merely on $\pi$.
A closer look, cf. the formulas in our paper \cite\aspheric,
shows that
the second cohomology group
$\roman H^2(\pi,\bold Z)$
admits the following description:
Write
$\zeta$
for the $\bold Z$-valued \lq\lq cocycle\rq\rq\
which assigns
$1$ to $[r]$
and
$0$ to the other generators of
$N\big /[F,N]$ and,
for $1 \leq j \leq n$,
write
$\zeta_j$
for the $\bold Z$-valued \lq\lq cocycle\rq\rq\
which assigns
$1$ to $[r_j]$
and
$0$ to the other generators of
$N\big /[F,N]$.
The
group
$\roman H^2(\pi,\bold Z)$
decomposes into a direct sum
of an infinite cyclic group, generated
by the class of $\zeta$,
and
$n$ finite cyclic groups
of orders respectively
$m_1,\dots, m_n$,
generated by the classes of
$\zeta_1,\dots,\zeta_n$, respectively.
In particular,
with the integers as coefficients,
unlike for the case of the fundamental group of a surface,
there is no canonical choice of
universal central extension.
On the other hand,
the second cohomology group
$\roman H^2(\pi,\bold R)$
is a one-dimensional real vector space
and hence there is a
{\it universal central extension\/}
$$
0
@>>>
\bold R
@>>>
\Gamma_{\bold R}
@>>>
\pi
@>>>
1
\tag3.1
$$
of $\pi$ by the reals
which does {\it not\/} depend on the choice of presentation.
In fact, the latter is determined by the requirement that
a  diagram of the kind
$$
\CD
0
@>>>
N\big / [F,N]
@>>>
\Gamma
@>>>
\pi
@>>>
1
\\
@.
@V{\lambda}VV
@VVV
@VV{\roman{Id}}V
@.
\\
0
@>>>
\bold R
@>>>
\Gamma_{\bold R}
@>>>
\pi
@>>>
1
\endCD
\tag3.2
$$
be commutative, where
$\lambda [r] = 1$ while,
for $j = 1,\dots,n$,
$\lambda [r_j]$
may be chosen arbitrarily;
whatever choice of the values
$\lambda [r_j]$, the resulting
extensions of $\pi$ by $\bold R$ will be {\it congruent\/}.
In other words, if the two extensions
$$
0
@>>>
\bold R
@>>>
\Gamma_1
@>>>
\pi
@>>>
1,
\quad
0
@>>>
\bold R
@>>>
\Gamma_2
@>>>
\pi
@>>>
1
$$
arise from homomorphisms $\lambda_1$ and $\lambda_2$, respectively,
as in (3.2) above, $\lambda_1$ and $\lambda_2$ determine a commutative
diagram
$$
\CD
0
@>>>
\bold R
@>>>
\Gamma_1
@>>>
\pi
@>>>
1
\\
@.
@V{\roman{Id}}VV
@VVV
@V{\roman{Id}}VV
@.
\\
0
@>>>
\bold R
@>>>
\Gamma_2
@>>>
\pi
@>>>
1
\endCD
$$
of group extensions.
Write $Z$ for the centre of $G$,
let $z$ be the Lie algebra of $Z$, and
let $X \in z$.
When $G$ is connected, $z$ coincides with the zentre of $g$
but in general $z$ equals the invariants for the induced
action of the group $\pi_0$ of components of $G$ on the centre of $g$.
Let $\roman{Hom}_X(\Gamma_{\bold R},G)$
denote the space of homomorphisms
$\phi$ from $\Gamma_{\bold R}$ to $G$
having the property that
$\phi (t[r]) = \roman{exp}(tX)$;
we assume $X$ chosen so that
$\roman{Hom}_X(\Gamma_{\bold R},G)$ is non-empty.
This assumption is topological in nature.
We briefly explain this below.
Let $\roman{Rep}_X(\Gamma_{\bold R},G)
=\roman{Hom}_X(\Gamma_{\bold R},G)\big / G$,
the resulting
{\it twisted representation space\/}.
The space
$\roman{Hom}_X(\Gamma_{\bold R},G)$
and hence
$\roman{Rep}_X(\Gamma_{\bold R},G)$
is unambigously defined,
independently of the choice
of presentation etc.,
since so is the central extension (3.1).
The space
$\roman{Rep}_X(\Gamma_{\bold R},G)$
is one of {\it projective\/} representations
of our  group $\pi$.
\smallskip
The
choice of generators
in $\Cal P$
identifies
the space $\roman{Hom}_X(\Gamma_{\bold R},G)$ with the pre-image
of
$\roman{exp}(X) \in G$, for the word map
$r$ from
$\roman{Hom}(F^\natural,G)$
to $G$
induced by the relator $r$, and we can play a similar
game as before, with the same choice
of $c \in C_2(F)$ so that $\partial c = r$
represents $\kappa \in \roman H_2(\pi,\bold R)$.
More precisely,
since
the centre of $G$ is contained in $O$,
the space $\roman{Hom}_X(\Gamma_{\bold R},G)$
arises as
the pre-image of $X \in z \subseteq O$
under the map
$\widehat r$ from $\Cal H^\natural(F,G)$ to $O$.
Furthermore,
in view of the Corollary to Lemma 1 in \cite\modus,
the map $\psi$ from $g$ to $g^*$
is regular at every point of the centre of $g$,
in fact, the restriction of $\psi$ to the centre
equals  the adjoint of
the given 2-form
whence
the space $\roman{Hom}_X(\Gamma_{\bold R},G)$
equals
the pre-image of
the adjoint
$X^{\sharp}\in g^*$
of $X$
under the momentum mapping $\mu$
from $\Cal M(\Cal P,G)$ to $g^*$.
Consequently the space $\roman{Rep}_X(\Gamma_{\bold R},G)$
is the corresponding reduced space,
for the coadjoint orbit
in $g^*$
consisting of the single point $X^{\sharp}$.
We already pointed out that,
when the standard homotopy operator
on forms on $g$ is taken,
the map $\psi$ from $g$ to
$g^*$ is in fact the adjoint of the chosen 2-form
$\cdot\,$ on $g$.
\smallskip
The ambiguity
with the choice of universal central extension over the integers
is of course resolved
by such a choice.
Extensions of this kind
have recently become of interest in the literature.
We therefore explain briefly the resulting representation
theory from our point of view:
An arbitrary central extension $\Gamma_{(b,\beta_1,\dots,\beta_n)}$
of $\pi$ by the integers
is
given by a presentation
$$
\langle x_1,y_1,  \dots,x_\ell,y_\ell,z_1,\dots,z_n, h;
[h,x_j],[h,y_j], [h,z_j], rh^{-b}, r_j h^{-\beta_j}
\rangle
$$
where
$$
r = \Pi [x_j,y_j] z_1 \dots z_n,
\  r_j = z_j^{m_j}
$$
as before
and where
the parameters
$b,\beta_1,\dots,\beta_n$
are arbitrary integers; they
correspond of course to the decomposition of $\roman H^2(\pi,\bold Z)$
mentioned earlier, and different choices
of these parameters may lead to the same group.
In particular, a group of the kind
$\Gamma_{(1,\beta_1,\dots,\beta_n)}$
fits into a central extension
$$
0
@>>>
\bold Z
@>>>
\Gamma_{(1,\beta_1,\dots,\beta_n)}
@>>>
\pi
@>>>
1
$$
in such a way that
the class $[r]$ of $r$ in $\Gamma_{(1,\beta_1,\dots,\beta_n)}$ is
identified with a generator of $\bold Z$.
Let $\roman{Hom}_X(\Gamma_{(1,\beta_1,\dots,\beta_n)},G)$
denote the space of homomorphisms
$\phi$ from $\Gamma_{(1,\beta_1,\dots,\beta_n)}$ to $G$
having the property that
$\phi ([r]) = \roman{exp}(X)$.
Let $\roman{Rep}_X(\Gamma_{(1,\beta_1,\dots,\beta_n)},G)
=\roman{Hom}_X(\Gamma_{(1,\beta_1,\dots,\beta_n)},G)\big / G$,
the resulting {\it twisted representation space\/}.
The choice of generators identifies
the space $\roman{Rep}_X(\Gamma_{\bold R},G)$
with the space
$\roman{Rep}_X(\Gamma_{(1,\beta_1,\dots,\beta_n)},G)$,
whatever
$(\beta_1,\dots,\beta_n)$.
\smallskip
When the parameters
$b,\beta_1,\dots,\beta_n, m_1,\dots, m_n$
satisfy certain numerical conditions,
$\Gamma_{(b,\beta_1,\dots,\beta_n)}$
is the fundamental group of a {\it Seifert fibre space\/}
which is  (i) closed (as a 3-manifold),
(ii) is an Eilenberg-Mac Lane space,
and has (iii) orientable decomposition surface,
cf. \cite\orvogzie.
By {\it symplectic reduction in finite dimensions\/},
we thus  obtain in particular spaces of representations
of fundamental groups of {\it all\/} Seifert fibre spaces
belonging to the class described above.
The significance of this remark has already been spelled out
in the Introduction.
\smallskip
In our situation,
the parameter $X$ is a certain topological characteristic class
of a principal $V$-bundle over the orbit space
$\Sigma$
with structure group $G$
associated with the representations
in
$\roman{Rep}_X(\Gamma_{\bold R},G)$,
cf. \cite\furustee.

\beginsection 4. The connected components

The reduced spaces arising from the above construction
will in general have more than one connected component.
We now explain briefly how these components arise
and how they can be labelled.
\smallskip
Let $C$ denote a finite cyclic group.
Then the space $\roman{Hom}(C,G)$
has finitely many connected components,
whence the above space
$
\roman{Hom}(F^\natural,G)
$
and therefore the various representation spaces
will have
finitely many
connected components,
certainly more than one
if $\pi$ has elliptic elements.
\smallskip
We describe some of the details in a special case:
Suppose $G$ compact and connected, and let $T$ be a maximal torus
in $G$, of rank say $r$.
Then the space
$\roman {Hom}(C,T)$
is a finite set,
consisting of $r |C|$ points.
Consequently
the space
$\roman {Hom}(C,G)$
will have
$r |C|$
connected components, each one of the form
$G/K$ for a closed subgroup $K$ of $G$.
The connected components, in turn,
correspond to conjugacy classes
of elements in $G$.
Thus a connected component of
$
\roman{Hom}(F^\natural,G)
$
and therefore of the various representation spaces
will be determined by a choice
of conjugacy class,
for each one of the generators
$z_1,\dots,z_n$ of finite order.
\smallskip
We now consider the special case of $G=U(n)$,
the unitary group.
A choice of
$X$
so that $\roman{Hom}_X(\Gamma_{\bold R},G)$
is non-empty
will correspond to a certain
holomorphic rank $n$
vector bundle $\zeta$ on $\Sigma$.
Let
$p_1,\dots,p_n$
be the distinct points
on $\Sigma$
arising as the images of the fixed points
of the elliptic transformations
$z_1,\dots,z_n$
in $\pi$.
A choice
of conjugacy class
for each $z_j$
corresponds to picking
a flag and a rational weight,
for each $p_j$;
in other words
a connected component
corresponds
to a parabolic structure
with rational weights
on $\zeta$. Our
twisted moduli space
$\roman{Rep}_X(\Gamma_{\bold R},G)$
will then
contain as top stratum
a homeomorphic
image of the stable part of
the corresponding {\smc Mehta-Seshadri}
\cite\mehtsesh\
moduli space of semistable
parabolic rank $n$ holomorphic vector
bundles
with rational weights
of degree determined by $X$.
It is likely that in fact
our construction
yields
all of these
moduli spaces,
not just the stable part.
We already pointed out that, for parabolic degree zero,
the spaces {\it are\/} known to be homeomorphic,
by a result of {\smc Mehta-Seshadri}, cf.
(4.1) and (4.3) in \cite\mehtsesh.
What is {\it new\/} here is that
important geometric information
about the {\smc Mehta-Seshadri}
moduli spaces
is obtained by {\it symplectic reduction, applied to
a smooth finite dimensional symplectic
manifold with a hamiltonian action of
the finite dimensional Lie group\/} $\roman U(n)$.

\beginsection 5. Applications

Suppose $G$ compact.
Recall that the notion of stratified symplectic space
has been introduced in \cite\sjamlerm.

\proclaim {Theorem 5.1}
With respect to the decomposition according to $G$-orbit types,
each connected component of the space
$\roman{Rep}(\pi,G)$
and, more generally,
each connected component of a
twisted representation space
$\roman{Rep}_X(\Gamma_{\bold R},G)$
inherits a structure of stratified symplectic space.
\endproclaim

In fact, the argument
for
the main result of
\cite\sjamlerm\
shows that
each connected
component of a reduced space of the kind considered
inherits a structure of stratified symplectic space.
In the setting
of \cite\sjamlerm\
the hypothesis of properness
is used {\it only \/}
to
guarantee that the reduced space is in fact connected.
In our situation, we know a priori
that the reduced space is connected.

\proclaim{Theorem 5.2}
Each stratum of (a connected component of) the space
$\roman{Rep}(\pi,G)$
and, more generally,
each stratum of (a connected component of) a
twisted representation space
$\roman{Rep}_X(\Gamma_{\bold R},G)$
has finite symplectic volume.
\endproclaim

The proof follows the same pattern
as that
for the argument for
(3.9) in \cite\sjamlerm.
There the unreduced symplectic manifold is assumed compact.
However the compactness of the zero level set suffices;
in our situation the zero level set {\it is\/} compact.
In fact, it suffices to prove the statement for the local model
in \cite\sjamlerm\
which looks like the reduced space of a unitary representation
of a compact Lie group, for the corresponding unique momentum
mapping having the value zero at the origin.
For the local model there is no difference between
(3.9) in \cite\sjamlerm\
and our situation.
Once the statement is established for the local model,
that of Theorem 5.2
follows since the reduced space
may be covered by finitely many open sets
having a local model of the kind described.
\smallskip
We mention two other consequences:

\proclaim{Corollary 5.3}
Each connected component
has a unique open, connected, and dense stratum.
\endproclaim

In fact, this follows at once from \cite\sjamlerm\ (5.9).
Likewise
\cite\sjamlerm\ (5.11) entails the following.

\proclaim{Corollary 5.4}
For each connected component,
the reduced Poisson algebra is symplectic, that is,
its only Casimir elements are the constants.
\endproclaim

\beginsection 6. Elliptic surfaces

Recall an elliptic surface
is a smooth compact complex surface $M$
with a proper surjective holomorphic map
$f \colon M \to \Sigma$
onto a complex curve $\Sigma$
such that the generic fibre is an elliptic curve.
Let $\{p_1,\dots,p_n\}\subseteq \Sigma$
be the finite set of non-regular values.
Then
$\Sigma$ may be viewed as an orbifold,
to each point $p_j$ the multiplicity of its
fibre being attached.
The map $f$ induces an isomorphism from
the fundamental group
$\pi = \pi_1(M)$ onto the orbifold fundamental
group of $(\Sigma,S)$.
Results of
{\smc Bauer} \cite\sbaueron\ relate
moduli spaces of
parabolic bundles
on $\Sigma$
of parabolic degree zero
with moduli spaces of
degree zero vector bundles
on $M$.
Bauer's construction
avoids {\smc Donaldson's} \cite\donalthr\
solution of the {\smc Kobayashi-Hitchin} conjecture.
In a follow up paper
\cite\einsther\
we shall study certain moduli spaces
over elliptic surfaces
from the symplectic point of view and in particular
offer a somewhat more intrinsic construction
of the moduli
space $\roman{Rep}_X(\Gamma_{\bold R},G)$
with its stratified symplectic structure,
viewed as a moduli space
of Einstein-Hermitian connections
on a certain projectively flat bundle
over $M$.
In particular, this will
give the symplectic counterpart of
Bauer's relationship
together with
(i) a compactification thereof and (ii)
an extension
 to arbitrary
parabolic degree.
A choice of
projective embedding of $M$
will then again correspond to a choice of universal central extension
of $\pi$
over
the integers.
Essentially different
choices of $\lambda$ in (3.2) will then correspond to
topologically inequivalent line bundles.
\bigskip

\centerline{References}
\medskip\noindent
References \cite\atiboton\  -- \cite\weinstwe\
are given in our paper \cite\modus.
\medskip

\widestnumber\key{999}
\ref \no \sbaueron
\by S. Bauer
\paper Parabolic bundles, elliptic surfaces,
and {\rm SU(2)}-representations of genus zero Fuchsian groups
\jour Math. Ann.
\vol 290
\yr 1991
\pages 509--526
\endref

\ref \no \baueroko
\by S. Bauer and C. Okonek
\paper The algebraic geometry of representation spaces
associated to Seifert fibred spaces
\jour Math. Ann.
\vol 286
\yr 1990
\pages 45--76
\endref

\ref \no \biereck
\by R. Bieri and B. Eckmann
\paper Groups with homological duality generalizing
Poincar\'e duality
\jour Inventiones Math.
\vol 20
\yr 1973
\pages 103-124
\endref

\ref \no \bodenone
\by H. U. Boden
\paper Unitary representations of Brieskorn spheres
\jour Duke Math. Journal
\vol 75
\yr 1994
\pages 193-220
\endref

\ref \no \kbrownon
\by K. S. Brown
\paper Groups of virtually finite dimension
\paperinfo
Proceedings of a symp. held at Durham in Sept. 1977, on
\lq \lq Homological and combinatorial techniques in group theory\rq \rq ,
ed. C.T.C. Wall,
London Math. Soc. Lecture Note Series 36
\publ
Cambridge Univ. Press
\publaddr Cambridge, U. K.
\yr 1979
\pages 27--70
\endref
\ref \no \donalthr
\by S. K. Donaldson
\paper Anti self-dual Yang-Mills connections over complex algebraic
surfaces and stable vector bundles
\jour Proc. Lond. Math. Soc.
\vol 50
\yr 1985
\pages 1--26
\endref

\ref \no \furustee
\by M. Furuta and B. Steer
\paper Seifert fibred homology 3-spheres and the Yang-Mills equations
on Riemann surfaces with marked points
\jour Adv. Math.
\vol 96
\yr 1992
\pages 38--102
\endref

\ref \no \guhujewe
\by K. Guruprasad, J. Huebschmann, L. Jeffrey, and A. Weinstein
\paper Groupoid representations and moduli spaces of parabolic bundles
\paperinfo in preparation
%\jour   \vol \yr \pages
\endref

\ref \no \modus
\by J. Huebschmann
\paper Symplectic and Poisson structures of certain moduli spaces
\paperinfo preprint, Nov. 1993, hep-th 9312112
\endref

\ref \no \aspheric
\by J. Huebschmann
\paper Cohomology theory of aspherical groups and of small cancellation
groups
\jour J. Pure and Applied Algebra
\vol 14
\pages 137--143
\yr 1979
\endref

\ref \no \einsther
\by J. Huebschmann
\paper Einstein bundles and extended moduli spaces
\paperinfo in preparation
\endref

\ref \no \huebjeff
\by J. Huebschmann and L. Jeffrey
\paper Group Cohomology Construction of Symplectic Forms
on Certain Moduli Spaces
\jour Int. Math. Research Notices
\vol 6
\yr 1994
\pages 245--249
\endref

\ref \no \kirkklas
\by P. Kirk and E. Klassen
\paper Representation spaces of Seifert fibered homology spheres
\jour Topology
\vol 30
\yr 1991
\pages 77--95
\endref

\ref \no \magnusbo
\by W. Magnus
\book Noneuclidean Tesselations And Their Groups
\publ Academic Press
\publaddr New York and London
\yr 1974
\endref

\ref \no \mehtrama
\by V. Mehta and A. Ramanathan
\paper Restriction of stable sheaves and representations
of the fundamental group
\jour Invent. Math.
\vol 77
\yr 1984
\pages 163--172
\endref

\ref \no \mehtsesh
\by V. Mehta and C. Seshadri
\paper Moduli of vector bundles on curves with parabolic structure
\jour Math. Ann.
\vol 248
\yr 1980
\pages 205--239
\endref

\ref \no \orvogzie
\by P. Orlik, E. Vogt, and H. Zieschang
\paper Zur Topologie gefaserter dreidimensionaler Mannigfaltigkeiten
\jour Topology
\vol 6
\yr 1967
\pages 49--64
\endref

\ref \no \patteone
\by S. J. Patterson
\paper On the cohomology of Fuchsian groups
\jour Glasgow Math. J.
\vol 16
\yr 1975
\pages 123--140
\endref

\ref \no \seshaone
\by C. S. Seshadri
\paper Spaces of unitary vector bundles on a compact Riemann surface
\jour Ann. of Math.
\vol 85
\yr 1967
\pages 303--336
\endref

\ref \no \seshaboo
\by C. S. Seshadri
\book Fibr\'es vectoriels sur les courbes alg\'ebriques
\bookinfo Ast\'erisque Vol. 96
\publ Soc. Math. de France
\yr 1982
\endref

\ref \no  \weiltwo
\by A. Weil
\paper Remarks on the cohomology of groups
\jour
Ann. of Math.
\vol 80
\yr 1964
\pages  149--157
\endref

\enddocument